\documentclass{amsart} 
\usepackage{amssymb}
\usepackage{amsmath}
\usepackage{amsfonts}

\begin{document}
\newtheorem{theo}{Theorem}[section]
\newtheorem{prop}[theo]{Proposition}
\newtheorem{lemma}[theo]{Lemma}
\newtheorem{exam}[theo]{Example}
\newtheorem{coro}[theo]{Corollary}
\theoremstyle{definition}
\newtheorem{defi}[theo]{Definition}
\newtheorem{rem}[theo]{Remark}


\newcommand{\Bb}{{\bf B}}
\newcommand{\Nb}{{\bf N}}
\newcommand{\Qb}{{\bf Q}}
\newcommand{\Rb}{{\bf R}}
\newcommand{\Zb}{{\bf Z}}
\newcommand{\Ac}{{\mathcal A}}
\newcommand{\Bc}{{\mathcal B}}
\newcommand{\Cc}{{\mathcal C}}
\newcommand{\Dc}{{\mathcal D}}
\newcommand{\Fc}{{\mathcal F}}
\newcommand{\Ic}{{\mathcal I}}
\newcommand{\Jc}{{\mathcal J}}
\newcommand{\Lc}{{\mathcal L}}
\newcommand{\Oc}{{\mathcal O}}
\newcommand{\Pc}{{\mathcal P}}
\newcommand{\Sc}{{\mathcal S}}
\newcommand{\Tc}{{\mathcal T}}
\newcommand{\Uc}{{\mathcal U}}
\newcommand{\Vc}{{\mathcal V}}

\newcommand{\ax}{{\rm ax}}
\newcommand{\Acc}{{\rm Acc}}
\newcommand{\Act}{{\rm Act}}
\newcommand{\ded}{{\rm ded}}
\newcommand{\Gm}{{$\Gamma_0$}}
\newcommand{\ID}{{${\rm ID}_1^i(\Oc)$}}
\newcommand{\PA}{{\rm PA}}
\newcommand{\ACA}{{${\rm ACA}^i$}}
\newcommand{\RefP}{{${\rm Ref}^*({\rm PA}(P))$}}
\newcommand{\RefS}{{${\rm Ref}^*({\rm S}(P))$}}
\newcommand{\Rfn}{{\rm Rfn}}
\newcommand{\tar}{{\rm Tarski}}
\newcommand{\UNFA}{{${\mathcal U}({\rm NFA})$}}

\author{Nik Weaver}

\title [Kinds of concepts]
       {Kinds of concepts}

\address {Department of Mathematics\\
          Washington University in Saint Louis\\
          Saint Louis, MO 63130}

\email {nweaver@math.wustl.edu}

\date{\em December 28, 2011}

\maketitle


The contrast between sets and proper classes is puzzling because our
naive notion of a collection seems to be a single notion, not two
separate notions. Mathematicians tend to be unclear about the exact
nature of this dichotomy. There obviously is a basic distinction
to be made, but one has trouble saying just what it is.

Some might say that proper classes are not genuine objects
at all, that talk about proper classes is merely a convenient abuse of
language that could, if one wanted, always be bypassed in favor of direct
reference to sets. This certainly eliminates the problem of having to
resolve the notion of a collection into two categories, but it does not
get to the heart of the matter because we then have to explain why some
concepts (like {\it natural number}) have extensions, i.e., determine
collections, but others (like {\it ordinal}) do not. There is clearly
something right about the idea that proper classes are not delimited
objects in the same way that sets are, but we need to understand what
this really means.

Another common response invokes the iterative conception of sets,
according to which sets are to be thought of as being formed in stages.
The idea is that there is no set of all ordinals, for example, because
there is no stage at which all the ordinals are available to be formed
into a set. The obvious problem for this explanation is how we are to
interpret the language about set formation. Surely it is not meant
literally, as an activity that can actually be performed. (Do sets exist
before they are formed? Could the same set be formed more than once?)
But how else could it be meant? This is a critical difficulty; if we
have no notion of what it means to form a set then saying that the
extension of the concept {\it ordinal} cannot be formed is not
illuminating. Again, there is clearly something right about the idea
that proper classes are somehow inexhaustible, but saying that
they are formed in stages does not give it proper expression.

The set/class phenomenon persists irrespective of one's views on the global
extent of mathematics. According to a finitist, individual natural numbers
exist but there is no set of all natural numbers; according to a predicativist,
individual real numbers exist but there is no set of all real numbers;
according to a set theoretic platonist, individual ordinals exist but there
is no set of all ordinals. In each case a concept that is accepted as
meaningful fails to aggregate, in some unclear sense.

\section{Indefinite extensibility}

The idea of a proper class seems to be related to Michael Dummett's notion
of indefinite extensibility. According to Dummett (\cite{D1}, p.\ 195),
\begin{quote}
a concept is indefinitely extensible if, for any definite characterisation
of it, there is a natural extension of this characterisation, which yields
a more inclusive concept; this extension will be made according to some
general principle for generating such extensions, and, typically, the
extended characterisation will be formulated by reference to the previous,
unextended, characterisation.
\end{quote}
Alternatively (\cite{D2}, p.\ 441),
\begin{quote}
an indefinitely extensible concept is one such that, if we can form a
definite conception of a totality all of whose members fall under that
concept, we can, by reference to that totality, characterize a larger
totality all of whose members fall under it.
\end{quote}
The concept {\it set} is considered
a typical example of indefinite extensibility
because any totality all of whose members are sets is itself a set, and so
if $x$ is the original totality then $x \cup \{x\}$ is a larger totality
all of whose members are sets. (Here we assume the axiom of foundation;
otherwise, replace ``set'' throughout the preceding comment with ``set
that is not an element of itself''.)

The passages quoted above appeal to our intuitive sense that proper classes
are in some way absolutely inexhaustible, and they are more cogent than the
iterative conception because what is being ``formed'' is a conception of a
set rather than the set itself. But their exact meaning is still
elusive. In the first one an indefinitely extensible concept is
contrasted with a ``definite characterisation of it'', which seems to
suggest that the original concept was ambiguous in some way. It sounds like
the process Dummett has in mind here is something like the following. Guided
by a vague, informal idea, we produce a sequence of precise formal
concepts that incompletely embody the original notion. None of them is
definitive, and indeed we have a general principle for converting any
given partial formalization into a broader, more inclusive formulation.
So there can be no ultimate precise, formal version of the original vague,
informal notion. The most we can achieve is an open-ended sequence of partial
formalizations.

This interpretation suggests a simple solution to the problem posed in
the introduction about why only some concepts have extensions. Perhaps
all we have to say is that concepts like {\it set} fail to have well-defined
extensions because they are vague?

However, this suggestion does not work. Call a concept {\it definite}
if every individual definitely either does or does not fall under it.
``Vague'' is conventionally taken to mean the opposite of this, i.e., as
implying that there could exist
individuals whose inclusion in the concept is undetermined. But that cannot
be what it means here because at least some concepts that certainly qualify
as indefinitely extensible are clearly not vague in this sense. For instance,
according to a finitist the concept {\it prime number} is indefinitely
extensible: given any finite set of prime numbers we have a finitary procedure
(multiply them together, add one, and factorize) which generates at least
one new prime number not in the set. Yet there is nothing ambiguous to a
finitist about what constitutes a prime number; indeed, we have a finitary,
mechanical procedure for testing primality. The assertion that any given
number is prime not only has a truth value, it has a decidable truth value.
So we have a concept which is both definite and indefinitely extensible.
(Probably this point is not controversial, but we need to be completely
clear about it because we want to leave open the possibility of reasoning
about concepts that are not definite.)

An analogous example can be given in the case of predicativism.
Predicativists accept constructions of countable length but reject
uncountable sets \cite{W1}.
So according to a predicativist the concept {\it irrational number} is
indefinitely extensible because any set of irrational numbers is countable
and hence can be diagonalized. But we can check with a computation of length
$2\omega$ whether a given infinite string of digits is eventually periodic, so
it is predicatively decidable whether such a string represents an irrational
number. Thus, we again have an indefinitely extensible concept whose
sense is understood with complete precision.
This example is slightly sharper than
the previous one because, while the finitist requires successively longer
computations to test for primality as the size of the candidate number
grows, the predicativist's procedure for testing irrationality always
involves a computation of the same length. So there is even less room to
argue that his grasp of the concept is changing in any way as his
repertoire of irrationals grows.

If there is any vagueness in cases like these, it resides not in the
application of the concept to any particular individual, but rather in
the question of where individuals falling under the concept are to be
sought. The finitist is, so to speak, not initially acquainted with all
the natural numbers and has no conception of a circumscribed arena in
which they appear. The predicativist is situated similarly with
respect to the real numbers, with the added feature that he not only
lacks initial familiarity with them, he does not even have a clear
generating procedure that would potentially produce all of them. Again,
the individuals falling under the concept are not to be found in any
circumscribed arena.

Now, contrary to what we suggested above, it is not obvious that concepts
which are ``vague'' only in the oblique sense of being uncircumscribed
cannot be said to have
well-defined extensions. But whether we want to say they have extensions
does not matter so much at this point because in any case
they do not have circumscribed extensions, and this ought to go a long
way toward explaining the set/class distinction. The problem is that
understanding what ``circumscribed'' means in this context does not seem
all that different from our original goal of understanding how sets differ
from classes. We already perceived that sets are limited
in some way that proper classes are not, and it is not clear that the
image of circumscription adds any precision to this thought.

The second passage quoted above frames the extensibility condition
in terms of our being able to enlarge any totality all of whose members
fall under the concept. This formulation squares better with the examples
just mentioned because there is less of an implication that the underlying
concept is evolving, but it also brings out more directly our need to
understand what it means to have a ``definite conception of a totality''.
In particular, we must ask how this differs from merely grasping a
well-defined concept. (Just how is having a definite conception of the
totality of prime numbers any different from knowing what prime numbers
are, or to put it differently, once we know what prime numbers are, what
more do we need in order to be able to say that we have a definite
conception of them as a totality? What is it that the finitist lacks
in this case?) This again seems
related to, perhaps even essentially identical to, our original question
about how sets and classes differ. The definition of indefinite
extensibility evidently presumes that this question has already been
answered. We now see that indefinite extensibility cannot be used as a
criterion for differentiating {\it set} and {\it class}; to the contrary,
we already need to be able to differentiate these concepts before we
can make sense of indefinite extensibility.

\section{Surveyable concepts}

The following idea might be helpful in this connection. Say that a
concept is {\it surveyable} if it is possible, in principle, to
exhaustively survey all of the individuals which fall under it. That
is, in principle we could perform to completion the task of examining,
one at a time, all the individuals falling under the concept. This is
morally equivalent to saying that the truth value of any first order
sentence whose variables range over the individuals falling under the
concept is, in principle, decidable, provided that all atomic formulas
appearing in the sentence are decidable. The idea is that if the range
of the variables is surveyable then we should be able to evaluate both
universal and existential quantifiers by direct inspection.

If a concept is surveyable, we will also say that its extension is
surveyable (leaving aside for now the question of whether there can
be unsurveyable extensions). We define a {\it set} to be the extension
of a surveyable concept.

The informal notion of surveyability is akin to the informal notions
of comptability and decidability; indeed, each of the three is in some
sense reducible to either of the others. However, computability
and decidability are usually considered against an (often unstated)
implicit assumption of finitism --- i.e., transfinite processes are
forbidden --- which we definitely do not want to insist on here. We
leave open the interpretation of the qualifier ``in principle'' in our
characterization of surveyability, noting only that it is generally accepted
that what is possible in principle goes well beyond what is possible in
practice, and that the conflict between finitism, predicativism, and
platonism can
be seen as stemming precisely from differing views on just how great this
discrepancy is, with the finitist accepting the possibility of tasks of any
finite length, but nothing more, and the predicativist accepting the
possibility of countably transfinite tasks, but nothing more.

It follows that we cannot fully cash out the notion of surveyability without
settling on one of these views over the others. But even before we do so,
it already seems more informative than ideas of ``circumscription''
and ``totality''. For instance, someone not previously familiar with the
paradoxes might find it difficult to gauge the meaning of the question
\begin{quote}
If we can decide whether any given individual falls under some concept,
does it follow that all the individuals falling under that concept
constitute a definite totality?
\end{quote}
and would probably be inclined to answer ``yes''. Phrased in terms of
surveyability, the analogous question is
\begin{quote}
If we can decide whether any given individual falls under some concept,
does it follow that we can exhaustively survey all the
individuals falling under that concept?
\end{quote}
and this is not only
more lucid, it also does a better job of locating the burden. The default
answer is clearly that being able to diagnose whether any given individual
has some property need not entail that we are able to inventory
all the individuals with that property; ergo, there
are, on the face of it, (at least) two distinct kinds of concepts: those
which are surveyable and those which are merely definite. This moves us
toward an explanation of the difference between sets and classes.

In particular, a resolution of the set theoretic paradoxes is available to
us now as a consequence of the premise that sets are surveyable collections,
not merely collections. The collections appearing in the standard
paradoxes (the collection of all ordinals, the collection of all sets,
etc.)\ are not truly paradoxical because they are not actually sets, and
they are not sets because they are not surveyable. If there is any
question as to what could render us unable, even in principle, to survey
all the sets there are, we just have to observe that being in the
position of having completed a survey of all the sets there are is not
a possible state of affairs. This is where the notion of indefinite
extensibility can be enlightening. No matter how many
sets we have managed to survey, the collection of all the sets we
have surveyed that do not contain themselves will necessarily
be a set that we have not yet surveyed. Thus, even in principle
there is no way we could survey all the sets there are.

So there are concepts which cannot possibly be surveyable, and thus
the dichotomy we proposed above between those concepts which are surveyable
and those which are not is significant. We may attribute the set theoretic
paradoxes to a mistaken implicit assumption that all concepts are surveyable;
once this assumption is denied these paradoxes evaporate.

Of course this is not a complete solution to the problem of the paradoxes
because it really just shifts the difficulty over to broader notions such
as {\it class} and {\it concept} which are not limited by the restriction
to surveyability. We will discuss these cases below. If our characterization
of sets as surveyable collections is granted, then what we have said so far
only defuses the paradoxes specifically involving sets.

Surveyability is a general conception that is compatible
with a range of views on the extent of mathematical possibility, so we
cannot expect it to lead to a resolution of the conflict between finitism,
predicativism, and platonism. Nonetheless, it may shed light on the
distinctions between these views. For instance, the fundamental question
separating platonism from predicativism is whether the surveyability of a
collection $x$ entails the surveyability of the concept
{\it subset of $x$}. This observation could at least help the platonist
understand the predicativist's reservations at accepting the power set of
$\omega$ as a valid set. One clearly has a naive intuition for the
surveyability of the natural numbers (specifically, a sense of their
sequential availability) that one does not have for $\Pc(\omega)$.

\section{Abstract objects}

Before going any further, we need to decide how we want to handle abstract
objects. It should be clear what we mean when we say that the prime numbers
less than 100 can be surveyed, and what is at stake when we ask whether all
the prime numbers can be surveyed. But this kind of language cannot be taken
literally, because numbers are supposed to be abstract objects that one
cannot directly inspect.

Of course there is no real problem in this instance because abstract numbers
are concretely represented by numerals. When we talk about testing a number
for primality, or generating a new prime number not in a given finite list,
we are thinking of working with numerals or something equivalently concrete,
not directly with abstract numbers (whatever those even are). This will be
our model for dealing with abstract objects generally: we will always demand
that we have some sort of concrete proxies for them. It is important to make
this point explicit because qualities like surveyability could conceivably
depend on the way we set up the system of proxies.

For our purposes, an individual counts as concrete if in principle it could
be directly examined. Thus, a finitist should accept any numeral as concrete,
and a predicativist should accept countable structures such as infinite
decimal expansions.

We will want to allow abstract objects to have more than one concrete
representation, so the appropriate conditions for an abstract concept to
be considered surveyable are that (1) we are able to survey the concrete
proxies for the abstract individuals falling under the concept and (2) we
are able to decide when two concrete proxies represent the same abstract
individual. For definite concepts, the corresponding requirements are that
(1) every concrete individual definitely is or is not a proxy and (2) any
two concrete proxies definitely do or do not represent the same abstract
individual.

The concept {\it set} illustrates our need to have this discussion because
it is not so obvious how we ought to represent sets concretely.
Representing them extensionally, by a list of their elements, is not
straightforwardly feasable if we want to account for sets of abstract
objects. In particular, this could be problematic for sets whose elements
are other sets. However, for the pure, well-founded sets of mathematics it
suggests the recursive procedure ``list the elements of the set, representing
each one by a vertex, then repeat this procedure for each vertex'', which
then leads us to represent these sets as rooted trees with no infinite paths
and no nontrivial automorphisms. In principle, trees can be
presented just as concretely as numerals can (a tree can be coded as an
order relation and then represented as an array of 0's and 1's), so they
satisfy our demand to have concrete proxies for sets. The criterion for
identity is the existence of a rooted tree isomorphism.

\section{Definite concepts}

We have said that a concept is definite if any individual definitely either
does or does not fall under it. Every surveyable concept is definite because
if we can survey all the individuals falling under some concept then we can
determine by inspection whether any given individual appears among them,
and so the assertion that it does appear must have a definite truth value.
But not every definite concept is surveyable.

Being definite but not surveyable is roughly the same as being indefinitely
extensible, but the correspondence is not exact. Using classical logic we can
deduce that any concept $C$ which is definite but not surveyable will have
the property that for any set $x$ all of whose members fall under $C$, there
is at least one individual falling under $C$ which does not belong to $x$.
However, in contrast to the indefinitely extensible case, we need not have
a ``natural extension'' or a ``general principle'' for enlarging sets
subordinate to $C$. The fact that any set subordinate to some
concept can be enlarged does not entail that this can necessarily be
accomplished by a uniform procedure.

On the other hand, if we are reasoning constructively then there is a
different discrepancy. To see this, suppose we could prove of some
definite concept $C$ that for every set $x$, not all the individuals
falling under $C$ belong to $x$. Then we would know that $C$ is not
surveyable, but we might not be able to prove it is indefinitely
extensible in Dummett's sense because intuitionistic logic does not let
us automatically infer that every set subordinate to $C$ is contained
in a larger set subordinate to $C$. To prove this we would have to be
able to generate new individuals falling under the concept, not merely
know that their nonexistence is absurd. We conclude that regardless of
which form of logic is in use, there could be definite concepts whose
extensions are not sets but which fail to be indefinitely extensible.

The question of which kind of reasoning is appropriate, classical or
constructive, may be contentious. One thing we can say here is that even
if we interpret the
logical constants constructively, classical logic will still be valid
for set theoretic sentences with bounded quantifiers (i.e., which only
quantify over surveyable collections). This is because, as we noted at
the beginning of this section, every instance of an atomic formula
asserting membership of an individual in a set is decidably true or
false, and as we noted at the beginning of Section 2, bounded
quantification preserves decidability. So set theoretic sentences with
bounded quantifiers have decidable truth values, and hence they
satisfy the law of excluded middle.

In contrast, the truth values of sentences whose variables range over
unsurveyable concepts generally cannot be settled by inspection, so if
the logical constants are interpreted constructively then excluded middle
is not assured in this setting. This brings us back to the question we
touched on earlier of whether every definite concept can be said to have
a well-defined extension; if not, then statements which quantify over
such concepts could not be assumed to have definite truth values and
constructive reasoning would be called for.

It is clear that if a concept is definite but not surveyable then we
cannot ``form a definite conception of a totality'' consisting of
precisely those individuals which fall under it. But that is just
because its extension is not a ``totality'', i.e., is not surveyable, not
because its extension is in any sense indefinite. Confusion can arise here
from equating ``totality'' with ``definite collection'', since this
implies that collections which are not totalities cannot be definite, and
hence that it must somehow be possible for definite concepts to have hazy
extensions. But once the distinction afforded by the notion of surveyability
is available there is no need to make this equation.

To the contrary, it is not clear what is wrong, if we are finitists,
with regarding the natural numbers as comprising a well-defined albeit
unsurveyable collection, or if we are predicativists or platonists, with
regarding the real numbers or the ordinals in the same way. Indeed, it is
not obvious what more could be reasonably demanded of a well-defined
collection, than that every object should definitely either belong or
not belong to it. Thus, we take the position that every definite concept
has a well-defined extension, and we define a {\it class} to be the
extension of a definite concept. We say it is {\it proper} if it is
not a set, i.e., it is not surveyable.

If we accept that classical logic is valid for statements which quantify over
a proper class, then we can affirm that the rooted tree representation of sets
discussed in the last section renders the concept {\it pure well-founded set}
definite. This is not immediate because testing whether a tree contains an
infinite path, or whether there exists an isomorphism between two trees,
requires us to quantify over all subsets of the vertex set (in the second
case, all subsets of the product of the two vertex sets). So we need
to determine whether classical logic holds when quantifying over all subsets
of a set. This certainly presents no problems for finitists or platonists,
because their views countenance power sets. But if
classical logic holds for statements which quantify over a definite
concept, then predicativists should not have any problems either. We can
model a subset of the set of vertices of a tree by coloring each vertex
either red or blue; if the original tree is surveyable then any purported
two-coloring of it can be vetted, so under this representation the concept
{\it subset of the set of vertices} is indeed definite. Thus, even for
a predicativist, assertions which quantify over all subsets of the set of
vertices of a tree have definite truth values, and so we conclude that the
concept {\it pure well-founded set} is definite.

The validity of classical logic also has implications for comprehension
axioms. Consider definitions of the form
$$z = \{y \in x: \phi(y)\}.$$
If $x$ is a set and $\phi$ is a set theoretic formula in which all
parameters are sets and all quantifiers range over sets, then $z$ will
also be a set; as we mentioned earlier in this section,
the truth value of $\phi(y)$ is decidable for all $y$, so it is possible
to survey $z$ by surveying $x$ and skipping over those individuals which
falsify $\phi$. Now consider the analogous claim for classes: that if $x$
is a class, all parameters in $\phi$ are classes, and all quantifiers in
$\phi$ range over classes, then $z$ will also be a class. Obviously this
holds in general only if classical logic is valid when quantifying over
definite concepts. The predicate $\phi$ represents a
definite concept precisely if it satisfies excluded middle.

\section{Indefinite concepts}

We resolved the set theoretic paradoxes by observing that the concept
{\it set} is not surveyable. The analogous paradoxes involving classes
would be similarly defeated if we found that the concept {\it class}
was not definite. This seems possible because the tree representation
that we used to show the concept {\it pure well-founded set} is definite
is not viable for proper classes (a point that would be completely missed
if we had ignored the need to concretely represent abstract objects).

At this point we have to ask whether it even makes sense to suppose that
a concept could fail to be definite. This could be seen as merely a
terminological issue which turns on how we define the word ``concept'',
but the substantive question that really matters is whether meaningful
atomic propositions which might not have a definite truth value
enter essentially into the global analysis of concepts. If they do, and
if we want to say that every meaningful predicate expresses a concept,
then we will have to accept that some concepts could be indefinite.

We claim that any global analysis of concepts must explicitly involve
the concept {\it valid proof}, and that this concept is indefinite.

Here we must distinguish between the syntactic notion of a valid proof
within a given formal system, which is not only definite but even decidable,
and the general semantic notion of a valid proof. One might hope to render
this distinction moot by capturing all semantically valid reasoning (or at
least all valid reasoning pertaining to the analysis of concepts) within a
single formal system. However, G\"odel's incompleteness theorem renders this
prospect doubtful. Any formal system that embodied only semantically valid
reasoning could presumably be recognized to do so, and could then, provided
it possessed sufficient number theoretic resources, be strengthened by adding
an arithmetical statement of its own consistency.

To elaborate on this last comment, we must indeed insist that any semantically
valid reasoning be, in principle, recognizably valid. Since the whole point
of deductive reasoning is to compel rational assent, it is essential that a
valid argument must be recognizable as such. However, this does not entail
that the concept {\it valid proof} is decidable; that would require that we
also be able to recognize any invalid reasoning to be invalid. The latter
property is not inherent in our concept of valid reasoning and cannot simply
be assumed.

To the contrary, the fact that one alleged proof can reference another
introduces a potential for circularity which could in some cases make
assigning a definite truth value to every assertion of the form ``$p$
is a valid proof'' problematic. For example, let $p$ be the
sentence ``By inspection, $q$ is a valid proof'' and let $q$ be
the sentence ``If $p$ were a valid proof then a falsehood would be
provable; therefore $p$ is not a valid proof.'' Since the validity of
$p$ entails the validity of $q$, which then entails the invalidity of $p$,
it is hard to avoid the conclusion that the reasoning exhibited in $q$ is
correct, i.e., that $q$ is a valid proof. But this shows that $p$ is valid,
which is absurd. Because of examples like these, we should
hesitate to assume that the concept {\it valid proof} is definite.

Although we are close to a paradox here, with careful analysis an outright
contradiction can be avoided. This analysis is presented in detail in
\cite{W3}; the point is that we can assume neither the general
validity of the law of excluded middle nor the general validity of the
inference of a statement $A$ from the premise that $A$ is provable. In short,
the latter inference hinges on the global soundness of all proofs and hence
becomes circular when it is itself used in a proof. As a result, there is a
subtle difference between asserting $A$ and asserting that $A$ is provable, and
in particular, asserting that $A$ is false, i.e., $A$ entails a contradiction,
is not the same as asserting that $A$ entails that a contradiction is
provable. This is enough to defuse the paradox in the example discussed
just above.

\section{Pseudo-concepts}

Next we must consider informal ideas like {\it true sentence} that cannot
be realized as concepts at all. It is easy to see that there can be no
predicate $T$, definite or indefinite, which globally satisfies Tarski's
biconditional $T(A) \leftrightarrow A$ for all sentences $A$: if there
were then we could construct a liar sentence which denies $T$ of itself,
leading immediately to a contradiction. Thus, our association of concepts
with predicates entails that {\it true sentence} cannot be a genuine concept.
We will call it a {\it pseudo-concept}.

The informal notion of truth is analyzed in \cite{W2}. We find that it
can be formalized in two distinct ways. The approach employed by
constructivists is to identify it with provability. Alternatively, we
can use Tarskian or Kripkean techniques to construct a predicate which
satisfies the condition $T(A) \leftrightarrow A$ for some limited family
of sentences $A$. This classical approach is subject to a ``revenge''
problem which shows that any such construction can always be broadened
to include more sentences within its scope. The situation here is aptly
described by Dummett's first characterization of indefinite extensibility
cited in Section 1, except that we prefer not to use the word ``concept''
for an informal idea that cannot be embodied in a well-defined predicate.
We might say that concepts like {\it pure well-founded set} have
indefinitely extensible reference, while pseudo-concepts like {\it true
sentence} have indefinitely extensible sense.

Thus, the pseudo-concept {\it true sentence} is classically realized as
a hierarchy of definite concepts. For each classical partial
truth predicate there is a corresponding liar sentence, but these are not
paradoxical; each liar sentence indeed does not fall under the predicate
to which it refers, but does fall under more inclusive truth predicates.

As we discussed in the last section, the constructive interpretation of
truth (viz., provability) evades the liar paradox by means of simultaneously
invalidating both the law of excluded middle and the inference of a statement
from the assertion that that statement is provable. However, there is still
an apparent difficulty associated with the classical hierarchy of partial
truth concepts, in the form of a ``partial liar'' sentence which asserts
that it does not fall under any concept in this hierarchy. This problem is
resolved by recognizing that the open-ended nature of our ability to construct
partial truth predicates reveals the concept {\it classical partial truth
predicate} to itself be constructive in nature. That is, we cannot say what
constitutes a classical partial truth predicate without referencing
provability. The precise definition we need is: a predicate $T$ is a partial
truth predicate if the implication $T(A) \to A$ is provable for any sentence
$A$. The appearance of provability in this formulation defeats the partial
liar paradox in the same way that the standard liar paradox is defeated when
we identify truth with provability. This appearance is unavoidable; there is
no way to define the notion of a partial truth predicate without it. We
cannot say something like ``$T$ is a partial truth predicate if the
implication $T(A) \to A$ holds for any sentence $A$'' because here ``holds''
is merely a synonym for ``is true''. Thus, we would need to already have a
global classical truth predicate in order to formulate a classical notion of
a partial truth predicate. We can only formulate this notion constructively.
Any particular classical partial truth predicate can still be defined without
invoking constructive ideas, but making global statements about all classical
partial truth predicates must be done constructively.

Thus, any global treatment, classical or constructive, of pseudo-concepts
like {\it true sentence} has to involve provability at some point.

\section{Falling under}

This conclusion has implications for the global analysis of concepts because
{\it falls under}, predicated of an object and a concept, is a
pseudo-relation in the same way that
{\it true sentence} is a pseudo-concept. Here the relevant biconditional is
$F(x,C) \leftrightarrow C(x)$ ($x$ falls under $C$ if and only if
$C(x)$), and it is now Russell's paradox which shows that there is no global
relation $F$ which satisfies the desired biconditional for all $x$ and $C$.

Like truth, falling under can be formalized in two ways, either constructively
in terms of $C$ provably holding of $x$ (i.e., the statement $C(x)$ is
provable), or classically, so as to satisfy the biconditional $F(x,C)
\leftrightarrow C(x)$ for some limited family of concepts $C$.

It also follows that informal notions such as {\it surveyable concept},
{\it definite concept}, {\it set}, and {\it class} whose formulation
involves falling under will themselves merely be pseudo-concepts, not
genuine concepts. This does not mean that these notions lie beyond formal
analysis. To the contrary, they can either be formalized constructively,
in terms of provability (e.g., a concept $C$ is constructively definite
if for any object $x$ the assertion $C(x) \vee \neg C(x)$ is provable) or
classically in some limited setting provided that we are able to define a
partial notion of falling under which satisfies the law
$F(x,C) \leftrightarrow C(x)$ for all concepts which appear
within that setting.

The class theoretic and concept theoretic paradoxes are now straightforwardly
resolved. If we are working classically within some limited setting, then
the paradoxical constructions will always take us out of that setting,
rendering the paradoxical reasoning invalid. Whereas if we formalize
falling under constructively then our inability to infer a statment $A$
from the assertion that $A$ is provable, together with the failure of the
law of excluded middle, is sufficient to disable the paradoxical reasoning.

In particular, the class theoretic paradoxes are not, as we speculated
earlier, to be resolved by observing that {\it class} is an indefinite
concept. They are resolved by observing that {\it class} is not a genuine
concept at all.


\end{document}